\documentclass[11pt,twoside,a4paper]{article}
\include{amsmath}
\usepackage{latexsym}
\usepackage{a4}
\usepackage{amstext}
\usepackage{amsfonts}\usepackage{amsmath}
\usepackage{amsthm}\usepackage{url}
\usepackage{amssymb}
\numberwithin{equation}{section}
\newtheorem{thm}{Theorem}[section]
\newtheorem{cor}[thm]{Corollary}

\theoremstyle{definition}
\newtheorem{defi}[thm]{Definition}

\begin{document}
\title{On the connection between Molchan-Golosov and Mandelbrot-Van\,Ness representations of fractional Brownian motion}
\author{C\'eline Jost\footnote{Tel.:\,+358-9-19151510; Fax:\,+358-9-19151400;
 E-mail address: celine.jost@helsinki.fi}}
\date{}\maketitle
\begin{small}Department of Mathematics and Statistics, P.O.\,Box 68 (Gustaf H\"allstr\"omin katu
2b), 00014 University of Helsinki, Finland.\end{small}
\begin{abstract}
We prove analytically a connection between the generalized Molchan-Golosov integral transform (see \cite{Jo}, Theorem 5.1) and the generalized Mandelbrot-Van\,Ness integral transform (see \cite{Pi2}, Theorem 1.1) of fractional Brownian motion (fBm). The former changes fBm of arbitrary Hurst index $K$ into fBm of index $H$ by integrating  over $[0,t]$,  whereas the latter requires integration over $(-\infty,t]$ for $t>0$. This completes an argument in \cite{Jo}, where the connection is mentioned without full proof.    
 \end{abstract}
\section{Introduction}
The \textit{fractional Brownian motion with Hurst index $H$} $\in(0,1)$, or \textit{$H$-fBm},  is the
continuous, centered Gaussian process $\left(B^H_t\right)_{t\in \mathbb{R}}$ with  $B^H_0 = 0$, a.s., and
\begin{equation*}\text{Cov}_{\mathbb{P}}\left(B^H_s,B^H_t\right)\ =\ \frac{1}{2}\left(|s|^{2H}\
+\ |t|^{2H}\ -\ |t-s|^{2H}\right),\ \ s,t\in\mathbb{R}.\label{covariance}\end{equation*} $H$-fBm is $H$-self-similar and has stationary increments. For $H=\frac{1}{2}$,  fractional Brownian motion is standard Brownian motion and denoted by $W$. FBm is interesting from a theoretical point of view, since it is fairly simple, but neither a Markov process, nor a semimartingale. Recently, the process has been studied extensively in connection to various applications, for example in finance and telecommunications. Important tools when working with fBm are its integral representations:
for a fixed Hurst index $K\in(0,1)$, on the one hand, there exists a $K$-fBm $\left(B^K_t\right)_{t\in\mathbb{R}}$, such that for all $t\in[0,\infty)$, we have that\begin{equation}B^H_t=C(K,H)\int^t_0(t-s)^{H-K}F\left(1-K-H,H-K,1+H-K,\frac{s-t}{s}\right)dB^{K}_s,\ a.s.,\label{1}\end{equation}
(see \cite{Jo}, Theorem 5.1). Here $C(K,H):=\frac{C(H)C(K)^{-1}}{\Gamma(H-K+1)}$, $C(H):=\Bigl(\frac{2H\Gamma\left(H+\frac{1}{2}\right)\Gamma\left(\frac{3}{2}-H\right)}{\Gamma(2-2H)}\Bigr)^{\frac{1}{2}}$, $\Gamma$ denotes the gamma function, and $F$ is Gauss' hypergeometric function. $B^K$ is unique (up to modification) on $[0,\infty)$.
On the other hand, there exists a unique (up to modification) $K$-fBm $\bigl(\tilde{B}^K_t\bigr)_{t\in\mathbb{R}}$, such that for all $t\in\mathbb{R}$, it holds that
\begin{equation}B^H_t\ =\ C(K,H)\int_{\mathbb{R}}\left((t-s)^{H-K}1_{(-\infty,t)}(s)-(-s)^{H-K}1_{(-\infty,0)}(s)\right)d\tilde{B}^K_s,\ a.s.\label{a1}\end{equation} 
(see \cite{Pi2}, Theorem 1.1). For $K=\frac{1}{2}$, (\ref{1}) corresponds to the Molchan-Golosov representation and (\ref{a1}) is the Mandelbrot-Van\,Ness representation of $H$-fBm (see \cite{Mo1} and \cite{Ma}, respectively).  The integrals in (\ref{1}) and (\ref{a1}) are fractional Wiener integrals. A priori, representations (\ref{1}) and (\ref{a1}) are very different. Indeed, the integrand in (\ref{1}) is a weighted fractional integral over $[0,t]$, whereas the integrand in (\ref{a1}) is a simple fractional integral over $\mathbb{R}$. Moreover, the filtrations generated by $\left(B^H_t\right)_{t\in[0,\infty)}$ and $\left(B^K_t\right)_{t\in[0,\infty)}$ coincide, but this is not the case for the natural filtrations of $\left(B^H_t\right)_{t\in\mathbb{R}}$ and $\bigl(\tilde{B}^K_t\bigr)_{t\in\mathbb{R}}$.\\
\\
In this work, we demonstrate how analytical facts of fractional integrals, combined with shifting properties of fBm, are used in order to establish a natural connection between the (generalized) Molchan-Golosov integral transform (\ref{1}) and the (generalized) Mandelbrot-Van\,Ness integral transform (\ref{a1}). More precisely, we show that the latter one emerges as a boundary case of a suitable time-shifted former one: Based on (\ref{1}), we construct a sequence of $H$-fBms which, for fixed $t$, converges in $L^2(\mathbb{P})$-sense  to (\ref{a1}). We will specify the rate of convergence. In particular, the generalized Mandelbrot-Van\,Ness representation is a consequence of the generalized Molchan-Golosov representation.\\
\\
The article is organized as follows. In Section 2, we first review the definition and some relevant facts of Gauss' hypergeometric function. Second, we define fractional integrals and derivatives over $\mathbb{R}$ and show the connection to fBm. Third, we recall the definition of the fractional Wiener integral over the real line. In Section 3, we derive the connection between the integral representations.  
\section{Preliminaries}
\subsection{Gauss' hypergeometric function}
The \textit{Gauss hypergeometric function} of parameters $a,b,c$ and variable $z\in\mathbb{R}$ is defined by the formal power series
 \begin{equation*}F(a,b,c,z)\ \ :=\  _2 F_1(a,b,c,z) \ :=\ \sum_{k=0}^{\infty}\frac{(a)_k (b)_k }{(c)_k }\frac{z^k}{k!},\label{e:2.12}\end{equation*}where $(a)_0:=1$ and $(a)_k:=a\cdot (a+1)\cdot\ldots\cdot (a+k-1)$, $k\in\mathbb{N}$.
We assume that $c\in\mathcal{A}:=\mathbb{R}\setminus\{\ldots,-2,-1,0\}$ for this to make sense. 
If $|z|<1$ or $|z|=1$ and $c-b-a>0$, then the series converges absolutely. If  furthermore $c>b>0$ for $z\in[-1,1)$ and $b>0$ for $z=1$, then it can be represented by the Euler integral 
 \begin{equation}F (a,b,c,z)\ =\  \frac{\Gamma(c)}{\Gamma(b)\Gamma(c-b)}\int^1_0 v^{b-1}(1-v)^{c-b-1}(1-zv)^{-a}dv\label{e:2.13}\end{equation} 
(see \cite{Er}, p. 59). If $c>b>0$ then the expression on the right-hand side of (\ref{e:2.13}) is  well-defined for all  $z\in(-\infty,1)$. For these parameters, we can hence extend the definition of $F$ to all  $z\in(-\infty,1)$ via (\ref{e:2.13}). In order to 
extend $F$  for fixed  $z\in(-\infty,1]$ to more general parameters, we consider Gauss' relations for neighbor functions.
Functions of type $F(a\pm m,b,c,z)$,  $F(a,b\pm m,c,z)$ and $F(a,b,c\pm m,z)$, $m\in\mathbb{N}$, are called \textit{contiguous} to $F(a,b,c,z)$. If $m=1$ then they are also called \textit{neighbors}. For any two  neighbors $F_1(z)$, $F_2(z)$ of  $F(a,b,c,z)$,  one has a linear relation of type \begin{equation}A(z) F(a,b,c,z)\ +\ A_1(z) F_1(z)\ +\ A_2(z) F_2(z)\ =\ 0,\label{e:2.14}\end{equation} where $A,A_1$ and $A_2$ are first degree polynomials with coefficients depending on $a,b$ and $c$.
 See \cite{Ab}, p. 558 for all $15$ relations. We use the neighbor relations in order to extend $F$ for $z\in(-\infty,1)$ to all parameters  such that $c\in\mathcal{A}$, and for $z=1$  to all parameters that satisfy $c,c-b-a\in\mathcal{A}$. Among the most important properties of $F$ are the symmetry  $$F(a,b,c,z)\ =\  F(b,a,c,z)$$ and  the reduction formula $$F(0,b,c,z)\ =\  F(a,b,c,0)\ =\ 1.$$ Also we have the linear transformation formula (see \cite{Ab}, p. 559)  
 \begin{equation}F(a,b,c,z)\ =\ (1-z)^{-a}F\left(a,c-b,c,\frac{z}{z-1}\right),\ z<1.\label{e:2.18}\end{equation}
In particular,  \begin{equation}F(a,b,b,z)\ =\ (1-z)^{-a},\ z<1.\label{e:2.181}\end{equation}  $F$ is smooth in $z$ and we have (see \cite{Ab}, p. 557)\begin{equation}\frac{d}{dz}F(a,b,c,z)\ =\ \frac{ab}{c}F(a+1,b+1,c+1,z)\label{missed}\end{equation} and 
\begin{equation}\frac{d}{dz}\bigl( z^{a} F(a,b,c,\pm z)\bigr)\ =\  a z^{a-1} F(a+1,b,c,\pm z).\label{e:2.20}\end{equation}
$F$ is left-continuous in $z=1$ and it holds that (see \cite{Er}, p. 9) \begin{equation}F(a,b,c,1)\ =\ \frac{\Gamma(c)\Gamma(c-b-a)}{\Gamma(c-b)\Gamma(c-a)}.\label{e:2.23}\end{equation}
Let $a,b>-1$ and $w<x<y$. Substituting $u:=(y-x)v+x$  in (\ref{e:2.13}) and using (\ref{e:2.18}) implies that\\
\\
$\int^y_x (y-u)^{b}(u-w)^c(u-x)^a du\  $ 
\begin{eqnarray}=\frac{\Gamma(a+1)\Gamma(b+1)}{\Gamma(a+2+b)}(y-x)^{1+a+b}(x-w)^c F\left(-c,a+1,a+2+b,\frac{y-x}{w-x}\right)\ \nonumber\\
=\frac{\Gamma(a+1)\Gamma(b+1)}{\Gamma(a+2+b)}(y-x)^{1+a+b}(y-w)^c F\left(-c,b+1,a+2+b,\frac{y-x}{y-w}\right). \label{e:2.263}\end{eqnarray}
If $x<y<w$, then we have that\\
\\
$\int^y_x (y-u)^{b}(w-u)^c(u-x)^a du\ $ \begin{eqnarray}=\frac{\Gamma(a+1)\Gamma(b+1)}{\Gamma(a+2+b)}(y-x)^{1+a+b}(w-x)^c F\left(-c,a+1,a+2+b,\frac{y-x}{w-x}\right) \label{e:2.264}\\
=\frac{\Gamma(a+1)\Gamma(b+1)}{\Gamma(a+2+b)}(y-x)^{1+a+b}(w-y)^c F\left(-c,b+1,a+2+b,\frac{y-x}{y-w}\right). \label{e:2.261}\end{eqnarray}
By linearly combining neighbor relations, we obtain relations of type (\ref{e:2.14}), where $F_1(z)$ and $F_2(z)$ are contiguous to $F(a,b,c,z)$ and $A$, $A_1$ and $A_2$ are polynomials of higher degree. An example for a contiguity relation is\\
\begin{equation}-c F(a,b-1,c,z)+(c-b+zb-za)F(a,b,c+1,z)+b(1-z)F(a,b+1,c+1,z)\  =\ 0.\label{igo}\end{equation}
It can be checked easily by using series.
 \subsection{Fractional calculus over the real line}
For more information on fractional calculus in the context of fractional Brownian motion, see \cite{Pi1} or \cite{Jo}. See \cite{Sa} for general information on fractional calculus.
\begin{defi}\label{d:3.1}
Let $\alpha>0$. The \textit{(right-sided) Riemann-Liouville fractional integral operator of order} $\alpha$  is defined by\begin{equation*}\bigl(\mathcal{I}^{\alpha}_{-}f\bigr)(s)\ :=\
 \frac{1}{\Gamma(\alpha)}\int^{\infty}_s f(u)(u-s)^{\alpha-1}du,\
 s\in\mathbb{R}.\label{e:3.2}\end{equation*} 
Let $\alpha\in(0,1)$. The \textit{(right-sided) Riemann-Liouville  fractional derivative operator of order} $\alpha$ is defined by $$\bigl(\mathcal{D}^{\alpha}_{-}f\bigr)(s) := \frac{-d}{ds}\bigl(\mathcal{I}^{1-\alpha}_{-}f\bigr)(s) = \frac{1}{\Gamma(1-\alpha)}\frac{-d}{ds}\int^{\infty}_{s}f(u)(u-s)^{-\alpha}du,\ s\in\mathbb{R}.$$
The \textit{(right-sided) Marchaud fractional derivative operator of order} $\alpha$   is defined by 
$$\left(\mathbf{D}^{\alpha}_{-}f\right)(s)\ :=\ \lim_{\epsilon\searrow 0}\left(\mathbf{D}^{\alpha}_{-,\epsilon}f\right)(s),\ a.e.\ s\in\mathbb{R},$$
where $$\left(\mathbf{D}^{\alpha}_{-,\epsilon}f\right)(s)\ :=\ \frac{\alpha}{\Gamma(1-\alpha)}\int^{\infty}_{\epsilon}\bigl(f(s)-f(u+s)\bigr)u^{-\alpha-1}du.$$Moreover,  \begin{equation*}\mathcal{D}^{0}_{-}f\ :=\ \mathbf{D}^0_{-}f\ :=\  f.\label{e:3.18}\end{equation*}
We set $$\mathcal{I}^{-\alpha}_{-}\ :=\ \mathbf{D}^{\alpha}_{-},\ \alpha\in[0,1).$$\end{defi}
If $\alpha\in(0,1)$,  $p\in[1,\frac{1}{\alpha})$ and $f\in L^p(\mathbb{R})$, then  $\mathcal{I}^{\alpha}_{-}f$ is well-defined (see \cite{Sa}, p. 94).  Clearly, $\mathcal{I}^1_{-}f$ exists for $f\in L^1(\mathbb{R})$. We have the composition formula
\begin{equation*}\mathcal{I}^{\alpha}_{-}\mathcal{I}^{\beta}_{-}f\ =\ \mathcal{I}^{\alpha+\beta}_{-}f ,\ f\in L^1(\mathbb{R}),\ \mathcal{I}^{\beta}_{-}f\in L^1(\mathbb{R}),\ \alpha,\beta\in(0,1].\label{semig}\end{equation*}
Furthermore,  \begin{equation} \mathcal{D}^{\alpha}_{-}\mathcal{I}^{\beta}_{-}f \ =\ \mathcal{I}^{\beta-\alpha}_{-}f, \ f\in L^1(\mathbb{R}),\ \mathcal{I}^{\beta}_{-}f\in L^1(\mathbb{R}),\ 0<\alpha\leq\beta\leq 1.\label{qwe1}\end{equation}
If  $f$ is piecewise differentiable with $\text{supp}(f')\subseteq (-\infty,y)$ for some $y\in\mathbb{R}$ and so that $\mathcal{D}^{\alpha}_{-}f$ exists, then 
$\mathbf{D}^{\alpha}_{-}f = \mathcal{D}^{\alpha}_{-}f$.\\
\\
Let $t\in\mathbb{R}\setminus\{0\}$ and set $1_{[0,t)}:=-1_{[t,0)}$ for $t<0$. For all $H\in(0,1)$, we have that $$\left(\mathcal{I}^{H-\frac{1}{2}}_{-}1_{[0,t)}\right)(s)\ =\ \frac{1}{\Gamma\left(H+\frac{1}{2}\right)}\left((t-s)^{H-\frac{1}{2}}1_{(-\infty,t)}(s)-(-s)^{H-\frac{1}{2}}1_{(-\infty,0)}(s)\right).$$
So for $K=\frac{1}{2}$, (\ref{a1}) can be written as  
\begin{equation}B^H_t\ =\ C(H)\int_{\mathbb{R}} \Bigl(\mathcal{I}^{H-\frac{1}{2}}_{-}1_{[0,t)}\Bigr)(s)d\tilde{W}_s,\ a.s.,\ t\in\mathbb{R}.\label{e:3.50}\end{equation}
\subsection{Fractional Wiener integrals over the real line}
We combine (\ref{e:3.50}) with the standard Wiener integral in order to obtain a meaning for the expression $\int_{\mathbb{R}}f(s)dB^H_s$ for suitable deterministic integrands $f$. For details on this topic, see \cite{Pi4}. For $H>\frac{1}{2}$, the space of integrands is given by
$$\Lambda(H)
 :=\ \biggl\{ f\in L^1(\mathbb{R})\ \bigg|\   \int_{\mathbb{R}} \left(\mathcal{I}_{-}^{H-\frac{1}{2}}f\right)(s)^2 ds\ <\ \infty\biggr\}.$$
For $H<\frac{1}{2}$, we have that 
$$\Lambda(H)
 :=\ \Bigl\{ f:\mathbb{R}\to\mathbb{R}\ \Big|\ \exists\ \phi_f\in L^2(\mathbb{R})\ \text{such that}\ f\ =\ \mathcal{I}^{\frac{1}{2}-H}_{-}\phi_f\Bigr\}.$$
For $f\in\Lambda(H)$, the \textit{(time domain) fractional Wiener integral} with respect to $B^H$  is defined by
$$I^H(f)\ :=\  \int_{\mathbb{R}} f(s)dB^H_s \ :=\ C(H)\int_{\mathbb{R}} \left(\mathcal{I}^{H-\frac{1}{2}}_{-} f\right)(s)d\tilde{W}_s.$$ 
\section{The connection between the integral transforms}
Consider transform  (\ref{1}). $B^K$ has stationary increments, so for every $s>0$, the  process \begin{equation*}B^{H,s}_t\ :=\ C(K,H)\int^t_0(t-u)^{H-K}\hat{F}\left(\frac{u-t}{u}\right)dB^{K}_{u-s},\  \ t\in[0,\infty),\label{3}\end{equation*} where  \begin{equation}\hat{F}(z)\ := \ F(1-K-H,H-K,1+H-K,z),\ \label{f}\end{equation}
is an $H$-fBm. The increments of $B^{H,s}$  are  stationary, hence the time-shifted process
\begin{equation*}Z^{H,s}_t \ :=\ B^{H,s}_{t+s}\ -\ B^{H,s}_s,\ t\in[-s,\infty),\label{ZHs=}\end{equation*} is an $H$-fBm.
By substituting $v:=u-s$, we obtain that
 $$Z^{H,s}_t\ =\  C(K,H)\left(\int^{t}_{-s}(t-v)^{H-K}\hat{F}\left(\frac{v-t}{v+s}\right)dB^K_v - \int^0_{-s}(-v)^{H-K}\hat{F}\left(\frac{v}{v+s}\right)dB^K_v\right),\ a.s.$$
 As $s\to\infty$, we formally obtain that
 \begin{equation*}Z^{H}_t\ :=\ Z^{H,\infty}_t\ :=\ C(K,H)\int_{\mathbb{R}}\left((t-v)^{H-K}1_{(-\infty,t)}(v)-(-v)^{H-K}1_{(-\infty,0)}(v)\right)dB^K_v\label{5}\end{equation*}
for $t\in\mathbb{R}$. Note that by definition, the processes $B^{H,s}$, $Z^{H,s}$ and $Z^{H}$ depend on $K$.
We can state and prove the following:
\begin{thm}\label{theorem1}For every $K\geq\frac{1}{2}$ and $t\in\mathbb{R}$, there exist constants $C_1=C_1(K,H,t)$ and $s_1=s_1(t)>0$, such that \begin{equation*}E\bigl[Z^{H,s}_t-Z^{H}_t\bigr]^2\ \leq\  C_1 s^{2H-2},\ s>s_1.\label{a}\end{equation*}
Moreover, for  every $K<\frac{1}{2}$ and $t\in\mathbb{R}$, there exist constants $C_2=C_2(K,H,t)$, $C_3=C_3(K,H,t)$ and $s_2=s_2(t)>0$ with \begin{equation*}E\bigl[Z^{H,s}_t-Z^{H}_t\bigr]^2\ \leq\ C_2 s^{2H-2}\ +\ C_3 s^{2K-2} ,\ s>s_2.\label{b}\end{equation*} \end{thm}
\begin{proof}Clearly, we can assume that $H\neq K$ and $t\neq 0$. Moreover, we assume that $t>0$. The result is derived similarly for $t<0$. Recall that $\hat{F}$ is defined in (\ref{f}), and  denote $$\Delta f^s_t(v)\ :=\ \left((t-v)^{H-K}\ -\  (-v)^{H-K}\right)1_{(-\infty,-s)}(v),$$
$$\Delta g^s_t(v)\ :=\ (t-v)^{H-K}\hat{F}\Bigl(\frac{v-t}{v+s}\Bigr)1_{(-s,t)}(v)\ -\ (-v)^{H-K}\hat{F}\Bigl(\frac{v}{v+s}\Bigr)1_{(-s,0)}(v),$$
$$\Delta h^s_t(v)\ :=\   (t-v)^{H-K}1_{(-s,t)}(v)\ -\ (-v)^{H-K}1_{(-s,0)}(v)$$
and  $\qquad\Delta k^s_t(v)\ :=\  k^s_t(v)\ -\ k^s_0(v)\ :=$
\begin{flushright}
$(t-v)^{H-K}\Bigl(\hat{F}\Bigl(\dfrac{v-t}{v+s}\Bigr)-1\Bigr)1_{(-s,t)}(v)\ -\ (-v)^{H-K}\Bigl(\hat{F}\Bigl(\dfrac{v}{v+s}\Bigr)-1\Bigr)1_{(-s,0)}(v).$\end{flushright}
For continuous $G$, set $ ^{\ast}G:=\max_{z\in[-1,0]}|G(z)|$ and $G^{\ast} := \max_{z\in[0,1]}|G(z)|$.\\
For $K=\frac{1}{2}$, we have by independence of increments of $B^K=W$, that
\begin{equation}\left(\frac{\Gamma\bigl(H+\frac{1}{2}\bigr)}{C(H)}\right)^2 E\bigl[Z^{H,s}_t-Z^H_t\bigr]^2\ =\  \int_{\mathbb{R}}\Delta f^s_t(v)^2 dv\ +\  \int_{\mathbb{R}}\Delta k_t^s(v)^2dv.\label{e:6.3}\end{equation}
Note that \begin{equation}\left(\sum_{i=1}^n a_i\right)^2\ \leq\ n\sum_{i=1}^n a_i^2,\ n\in\mathbb{N}.\label{sumequation}\end{equation}
So for $K\neq\frac{1}{2}$, we obtain by using (\ref{sumequation}) with $n=2$, that  \begin{eqnarray}\frac{1}{2C(K,H)^2}\cdot E\bigl[Z^{H,s}_t-Z^H_t\bigl]^2\ \leq\  E\left[\int_{\mathbb{R}}\Delta f^s_t(v)dB^K_v\right]^2 +  E\left[\int_{\mathbb{R}}\Delta k^s_t(v) dB^K_v\right]^2\nonumber\\ 
=C(K)^2 \left(\int_{\mathbb{R}}\left(\mathcal{I}^{K-\frac{1}{2}}_{-} \Delta f_t^s\right)(v)^2 dv +  \int_{\mathbb{R}}\left(\mathcal{I}^{K-\frac{1}{2}}_{-}\Delta k_t^s\right)(v)^2 dv\right).\label{e:6.4}\end{eqnarray}
\textbf{1.} For all $K\in(0,1)$, we show that there exists a constant $c_1(K,H,t)$, such that \begin{equation}\int_{\mathbb{R}}\left(\mathcal{I}^{K-\frac{1}{2}}_{-} \Delta f_t^s\right)(v)^2 dv\  \leq\ c_1(K,H,t) s^{2H-2},\ s>t.\label{e:6.20}\end{equation}  
For all $K\in(0,1)$, we  have that
$$\left(\mathcal{I}^{K-\frac{1}{2}}_{-}\Delta f_t^s\right)(v)
 =\ \frac{(-s-v)^{H-\frac{1}{2}}}{\Gamma\bigl(K+\frac{1}{2}\bigr)}\Bigl(G_0\Bigl(\frac{-s-v}{t-v}\Bigr)\ -\ G_0\Bigl(\frac{-s-v}{-v}\Bigr)\Bigr)1_{(-\infty,-s)}(v),$$
where $$G_0(z)\ :=\ z^{K-H}F\left(K-H,K-\frac{1}{2},K+\frac{1}{2},z\right).$$ 
For $K=\frac{1}{2}$, this is trivial. For $K>\frac{1}{2}$, this follows from (\ref{e:2.264}). For $K<\frac{1}{2}$, it follows from the fact that $\left(\mathcal{I}^{K-\frac{1}{2}}_{-}\Delta f_t^s\right)(v)=\frac{-d}{dv}\left(\mathcal{I}^{K+\frac{1}{2}}_{-}\Delta f_t^s\right)(v)$ by using (\ref{e:2.264}), (\ref{e:2.20}) and (\ref{igo}).
By using (\ref{e:2.20}), we have that
\begin{eqnarray}\frac{dG_0}{dz}(z) &=&  (K-H)z^{K-H-1}F\left(K-H+1,K-\frac{1}{2},K+\frac{1}{2},z\right)\nonumber\\  &=:&  (K-H)z^{K-H-1}G_1(z).\label{ucdc}\end{eqnarray} By the mean value theorem, there exists $\theta_{v,s}\in\bigl(\frac{-s-v}{t-v},\frac{-s-v}{-v}\bigr)\subset (0,1)$, such that\\
\\
$\left(\mathcal{I}^{K-\frac{1}{2}}_{-}\Delta f_t^s\right)(v)^2\ =\ \Bigl(\frac{K-H}{\Gamma\bigl(K+\frac{1}{2}\bigr)}\Bigr)^2 (-s-v)^{2H-1}\left(\frac{t(-v-s)}{-v(t-v)}\right)^2 \theta_{v,s}^{2(K-H-1)}G_1^2(\theta_{v,s})1_{(-\infty,-s)}(v)\ \leq$\\
\\
$\Bigl(\frac{K-H}{\Gamma\left(K+\frac{1}{2}\right)}\Bigr)^2 G_1^{\ast 2}t^2 \max\left(1,(\frac{t+s}{s})^{2(H-K)}\right)(-s-v)^{2K-1}(-v)^{2(H-K-1)}\cdot 1_{(-\infty,-s)}(v)$.\\
\\
Hence, from (\ref{e:2.261}) and (\ref{e:2.18}), we obtain that $\int_{\mathbb{R}}\left(\mathcal{I}^{K-\frac{1}{2}}_{-}\Delta f_t^s\right)(v)^2 dv\ \leq$\\
\\
$\Bigl(\frac{K-H}{\Gamma\left(K+\frac{1}{2}\right)}\Bigr)^2 G_1^{\ast 2}t^2 \max\left(1,(\frac{t+s}{s})^{2(H-K)}\right)\cdot $\vspace{-3mm}\begin{flushright}$\lim_{x\to-\infty}\frac{1}{2K}(-x-s)^{2K}s^{2(H-K-1)}\left(\frac{-x}{s}\right)^{-2K}F\left(2H-1,2K,2K+1,\frac{x+s}{x}\right)\ =$\end{flushright}\vspace{-2mm}
$\Bigl(\frac{K-H}{\Gamma\left(K+\frac{1}{2}\right)}\Bigr)^2 G_1^{\ast 2}\frac{t^2}{2K}F(2H-1,2K,2K+1,1)\max\left(1,(\frac{t+s}{s})^{2(H-K)}\right)s^{2H-2}$.\\
\\
By using (\ref{e:2.23}), we obtain (\ref{e:6.20}) with $$c_1(K,H,t)\ =\ \left(\frac{K-H}{\Gamma(K+\frac{1}{2})}\right)^2 \frac{G_1^{\ast 2}4\Gamma(2K)\Gamma(2-2H)}{\Gamma(2K-2H+2)}t^2.$$
\textbf{2.} Fix $d>0$. For $K\geq\frac{1}{2}$, we show that there exists a constant $c_2(K,H,t,d)$, such that  \begin{equation}\int_{\mathbb{R}}\left(\mathcal{I}^{K-\frac{1}{2}}_{-}\Delta k^s_t\right)(v)^2 dv\ \leq\ c_2(K,H,t,d)s^{2H-2},\ s>2t+4d+1.\label{e:6.21}\end{equation}
For $K<\frac{1}{2}$, we show that there exist constants $c_3(K,H,t,d)$ and $c_4(K,H,t,d)$ with  \begin{equation}\int_{\mathbb{R}}\left(\mathcal{I}^{K-\frac{1}{2}}_{-}\Delta k^s_t\right)(v)^2 dv \leq c_3(K,H,t,d)s^{2H-2} + c_4(K,H,t,d)s^{2K-2},\ s>2t+4d+1.\label{e:6.22}\end{equation} Let $s>2d$. We have that
 \begin{eqnarray*}\Delta k^s_t(v)\ &=&\ k^s_t(v) 1_{[-d,t)}(v)\ -\ k^s_0(v) 1_{[-d,0)}(v)\ +\ \Delta k^s_t(v)1_{[\frac{-s}{2},-d)}(v)\\
 & & \ +\ \ \Delta g^s_t(v)1_{(-s,\frac{-s}{2})}(v)\ -\ \Delta h^s_t(v)1_{(-s,\frac{-s}{2})}(v).\end{eqnarray*}
First, let $K=\frac{1}{2}$. Note that $(-s,\frac{-s}{2})\cap [\frac{-s}{2},-d)\cap [-d,t)=\emptyset$. By using (\ref{sumequation}) with $n=2$, we have that   \begin{eqnarray}\frac{1}{2}\int_{\mathbb{R}}\Delta k^s_t(v)^2 dv &\leq&  \int^t_{-d}k^s_t(v)^2 dv\ +\ \int^0_{-d}k^s_0(v)^2 dv \ +\ \int^{-d}_{\frac{-s}{2}}\Delta k^s_t(v)^2 dv\nonumber\\
& &\ +\ \int^{\frac{-s}{2}}_{-s}\Delta g^s_t(v)^2 dv\ +\ \int^{\frac{-s}{2}}_{-s}\Delta h^s_t(v)^2 dv. \label{e:6.5}
\end{eqnarray}
Second, let $K>\frac{1}{2}$. Then $\Gamma\left(K-\frac{1}{2}\right)\left(\mathcal{I}^{K-\frac{1}{2}}_{-}\Delta k^s_t\right)(v)\ =$\\
\\
$\int^t_{-d}k^s_t(u)(u-v)^{K-\frac{3}{2}}du \cdot 1_{(-\infty,-d)}(v)\ +\ \int^t_{v}k^s_t(u)(u-v)^{K-\frac{3}{2}}du \cdot 1_{[-d,t)}(v)$\\
\\
$ -\ \int^0_{-d}k^s_0(u)(u-v)^{K-\frac{3}{2}}du \cdot 1_{(-\infty,-d)}(v)\ -\ \int^0_{v}k^s_0(u)(u-v)^{K-\frac{3}{2}}du\cdot 1_{[-d,0)}(v)\ $\\
\\ 
$+ \int^{-d}_{\frac{-s}{2}}\Delta k^s_t(u)(u-v)^{K-\frac{3}{2}}du \cdot 1_{(-\infty,\frac{-s}{2})}(v) + \int^{-d}_{v}\Delta k^s_t(u)(u-v)^{K-\frac{3}{2}}du \cdot 1_{[\frac{-s}{2},-d)}(v)\ $\\
\\
$+ \int^{\frac{-s}{2}}_{-s}\Delta g^s_t(u)(u-v)^{K-\frac{3}{2}}du \cdot 1_{(-\infty,-s)}(v) + \int^{\frac{-s}{2}}_{v}\Delta g^s_t(u)(u-v)^{K-\frac{3}{2}}du\cdot 1_{[-s,\frac{-s}{2})}(v)\ \ $\\
\\
$- \int^{\frac{-s}{2}}_{-s}\Delta h^s_t(u)(u-v)^{K-\frac{3}{2}}du \cdot 1_{(-\infty,-s)}(v) -\int^{\frac{-s}{2}}_{v}\Delta h^s_t(u)(u-v)^{K-\frac{3}{2}}du \cdot 1_{[-s,\frac{-s}{2})}(v)$\\
\\
$=:\ A_1(v) + B_1(v) + A_2(v) + B_2(v) + A_3(v) +B_3(v) + A_4(v) + B_4(v) + A_5(v) +B_5(v)$.\\
\\
Hence, by using (\ref{sumequation}) with $n=5$, we obtain that
\begin{equation}\frac{\Gamma\left(K-\frac{1}{2}\right)^2}{5}\int_{\mathbb{R}}\left(\mathcal{I}^{K-\frac{1}{2}}_{-}\Delta k^s_t\right)(v)^2 dv \ \leq\   \sum_{i=1}^{5}\int_{\mathbb{R}}A_i(v)^2 dv\ +\ \sum_{i=1}^5\int_{\mathbb{R}} B_i(v)^2 dv.\label{e:6.6}\end{equation}
Third, let $K<\frac{1}{2}$. Let  $x\leq v<y$ and   $f$ be  differentiable with $f'\in L^1[x,y]$. From (\ref{qwe1}), it follows that $\left(\mathcal{I}^{K+\frac{1}{2}}_{-}\mathcal{D}^{1}_{-}f1_{(x,y)}\right)(v)\ =\ \left(\mathcal{D}^{\frac{1}{2}-K}_{-}\mathcal{I}^1_{-}\mathcal{D}^{1}_{-}f1_{(x,y)}\right)(v)\ =$\vspace{-3mm}
\begin{flushright}$
\left(\mathcal{D}^{\frac{1}{2}-K}_{-} f1_{(x,y)}\ -\ f(y)1_{(x,y)}\right)(v)\ =\
  \left(\mathcal{D}^{\frac{1}{2}-K}_{-}f1_{(x,y)}\right)(v)\ -\  \frac{f(y)}{\Gamma\left(K+\frac{1}{2}\right)}(y-v)^{K-\frac{1}{2}}$.\end{flushright}
Since $\lim_{v\nearrow t}k^s_t(v)=\lim_{v \nearrow 0}k^s_{0}(v)=0$, we obtain that\\
\\
 $\Gamma\left(K-\frac{1}{2}\right)\left(\mathbf{D}^{\frac{1}{2}-K}_-\Delta k^s_t\right)(v)\ =\ \Gamma\left(K-\frac{1}{2}\right)\left(\mathcal{D}^{\frac{1}{2}-K}_{-}k^s_t 1_{(-d,t)}\right)(v)\ -$\\
\\
$\Gamma\left(K-\frac{1}{2}\right)\left(\mathcal{D}^{\frac{1}{2}-K}_{-}k^s_0 1_{(-d,0)}\right)(v)\ +\  \Gamma\left(K-\frac{1}{2}\right)\left(\mathcal{D}^{\frac{1}{2}-K}_{-}\Delta k^s_t 1_{(\frac{-s}{2},-d)}\right)(v)\ +$\\
\\
$\Gamma\left(K-\frac{1}{2}\right)\left(\mathcal{D}^{\frac{1}{2}-K}_{-}\Delta g^s_t 1_{(-s,\frac{-s}{2})}\right)(v)\ -\ \Gamma\left(K-\frac{1}{2}\right)\left(\mathcal{D}^{\frac{1}{2}-K}_{-}\Delta h^s_t 1_{(-s,\frac{-s}{2})}\right)(v)\ =$\\
\\
$\int^t_{-d}k^s_t(u)(u-v)^{K-\frac{3}{2}}du\cdot  1_{(-\infty,-d)}(v)\ +\ \Gamma\left(K-\frac{1}{2}\right)\left(\mathcal{I}^{K+\frac{1}{2}}_{-}\frac{-d}{du}k^s_t 1_{(-d,t)}\right)(v)\cdot 1_{[-d,t)}(v)$\\
\\
$-\ \int^0_{-d}k^s_0(u)(u-v)^{K-\frac{3}{2}}du\cdot 1_{(-\infty,-d)}(v)\ -\ \Gamma\left(K-\frac{1}{2}\right)\left(\mathcal{I}^{K+\frac{1}{2}}_{-}\frac{-d}{du}k^s_0 1_{(-d,0)}\right)(v)\cdot 1_{[-d,0)}(v)$\\
\\
$+\ \int^{-d}_{\frac{-s}{2}}\Delta k^s_t(u)(u-v)^{K-\frac{3}{2}}du\cdot 1_{(-\infty,\frac{-s}{2})}(v)\ +$  \vspace{-3mm}\begin{flushright}$ \Gamma\left(K-\frac{1}{2}\right)\left(\mathcal{I}^{K+\frac{1}{2}}_{-}\frac{-d}{du}\Delta k^s_t 1_{(\frac{-s}{2},-d)}\right)(v)\cdot 1_{[\frac{-s}{2},-d)}(v)$ \end{flushright}\vspace{-3mm}
$+\ \int^{\frac{-s}{2}}_{-s}\Delta g^s_t(u)(u-v)^{K-\frac{3}{2}}du \cdot 1_{(-\infty,-s)}(v)\ +$\vspace{-3mm} \begin{flushright}$\Gamma\left(K-\frac{1}{2}\right)\left(\mathcal{I}^{K+\frac{1}{2}}_{-}\frac{-d}{du}\Delta g^s_t 1_{(-s,\frac{-s}{2})}\right)(v)\cdot 1_{[-s,\frac{-s}{2})}(v)$ \end{flushright}\vspace{-3mm}
$-\  \int^{\frac{-s}{2}}_{-s}\Delta h^s_t(u)(u-v)^{K-\frac{3}{2}}du \cdot 1_{(-\infty,-s)}(v)\ -$  \vspace{-3mm}\begin{flushright}$\Gamma\left(K-\frac{1}{2}\right)\left(\mathcal{I}^{K+\frac{1}{2}}_{-}\frac{-d}{du}\Delta h^s_t 1_{(-s,\frac{-s}{2})}\right)(v)\cdot 1_{[-s,\frac{-s}{2})}(v)$ \end{flushright}\vspace{-3mm}
$+\  \dfrac{\Delta k^s_t(-d)}{K-\frac{1}{2}}(-d-v)^{K-\frac{1}{2}}\cdot 1_{[\frac{-s}{2},-d)}(v)
\ +\ \dfrac{\Delta k^s_t\left(\frac{-s}{2}\right)}{K-\frac{1}{2}}\left(\frac{-s}{2}-v\right)^{K-\frac{1}{2}}\cdot 1_{[-s,\frac{-s}{2})}(v)$\\
\\
$=:\ A_1(v)+C_1(v)+A_2(v)+C_2(v)+A_3(v)+ C_3(v) + A_4(v) + C_4(v) + A_5(v) + C_5(v) + D(v) + E(v)$.\\
\\
Hence, (\ref{sumequation}) with $n=6$ yields
\begin{eqnarray}\frac{\Gamma\left(K-\frac{1}{2}\right)^2}{6}\int_{\mathbb{R}}\left(\mathcal{I}^{K-\frac{1}{2}}_{-}\Delta k^s_t\right)(v)^2 dv &\leq&  \sum_{i=1}^{5}\int_{\mathbb{R}}A_i(v)^2 dv\ +\ \sum_{i=1}^5 \int_{\mathbb{R}}C_i(v)^2 dv\nonumber\\
&&+\ \int_{\mathbb{R}}D(v)^2 dv\ +\ \int_{\mathbb{R}}E(v)^2 dv.\label{e:6.7}\end{eqnarray}
Next we estimate the integrals on the right-hand sides of (\ref{e:6.5}), (\ref{e:6.6})  and (\ref{e:6.7}). In what follows, $B$ denotes the beta function.\\
\\
1. Estimation of  $\int_{-d}^{t}k^s_{t}(v)^2 dv$,  $\int_{\mathbb{R}}A_{1}(v)^2 dv$,  $\int_{\mathbb{R}}B_1(v)^2 dv$ and  $\int_{\mathbb{R}}C_1(v)^2 dv$.\\
Let $u\in(-d,t)$. By the mean value theorem, there exists $\theta_{s,u}\in(\frac{u-t}{u+s},0)\subset\ (\frac{-d-t}{-d+s},0)$, such that\\
\\
$|k^s_t(u)|=(t-u)^{H-K}\bigl(\frac{t-u}{u+s}\bigr)\Big|\frac{d\hat{F}}{dz}(\theta_{s,u})\Big|\leq 
(t-u)^{H-K+1}(s-d)^{-1}\max_{z\in(\frac{-d-t}{-d+s},0)}\Big|\frac{d \hat{F}(z)}{dz}\Big|.$\\
\\
We assume that $s>2d+t$. Then, it follows from (\ref{missed}), that \begin{equation*}|k^s_t(u)|\ \leq\ (t-u)^{H-K+1}(s-d)^{-1}  \frac{|1-K-H||H-K| ^{\ast}G_2}{|1+H-K|},\label{proof1eq1}\end{equation*} where
 $$G_2(z)\ :=\ F(2-K-H,H-K+1,2+H-K,z).$$
For $K=\frac{1}{2}$, we obtain that $\int^t_{-d}k^s_t(u)^2 du\ \leq\ \frac{ ^\ast\! G^2_2  (t+d)^{2H+2}}{2}s^{-2}$.\\
\\
Denote $$G_3(z)\ :=\ F\left(\frac{3}{2}-K,H-K+2,H-K+3,z\right).$$
 For $K\neq \frac{1}{2}$, we obtain by using (\ref{e:2.263}), that
$\int_{\mathbb{R}}A_1(v)^2 dv\ \leq \frac{2 ^{\ast}\!G^2_2 G^{\ast 2}_3(t+d)^{2H+2}}{(1+H-K)^2 (1-K)}s^{-2}$.\\
\\
Also, for $K>\frac{1}{2}$, we have that $\int_{\mathbb{R}}B_1(v)^2 dv \ \leq\ \frac{2 ^{\ast}\!G^2_2 B\bigl(H-K+2,K-\frac{1}{2}\bigr)^2}{(1+H-K)^2}(t+d)^{2H+2}s^{-2}$.\\
\\
Furthermore, it holds that \\
\\
 $\big|\frac{d}{du}k^s_t(u)\big|\ \leq\ |K-H|(t-u)^{H-K-1}\bigl(\frac{t-u}{u+s}\bigr)\Big|\frac{d\hat{F}}{dz}(\theta_{s,u})\Big|\ +\ (t-u)^{H-K}\Big|\frac{d\hat{F}}{dz}\bigl(\frac{u-t}{u+s}\bigr)\Big|\frac{s+t}{(u+s)^2}\ \leq$\\
\\
$\left(\frac{|1-K-H||H-K| ^{\ast}\!G_2}{|1+H-K|}\right) \left(\frac{|K-H|}{s-d}\ +\ \frac{s+t}{(s-d)^{2}}\right)(t-u)^{H-K}.$\\
\\
Hence, for $K<\frac{1}{2}$, we obtain by using (\ref{sumequation}) with $n=2$, that \\
\\
$\int_{\mathbb{R}}C_1(v)^2 dv\ \leq \frac{ 68 \cdot ^{\ast}\! G^2_2  B\left(H-K+1,K+\frac{1}{2}\right)^2}{(1+H-K)^2 \left(K-\frac{1}{2}\right)^2}(t+d)^{2H+2}s^{-2}$.\\
\\
2. Estimation of  $\int^0_{-d}k^s_{0}(v)^2 dv$,  $\int_{\mathbb{R}}A_{2}(v)^2 dv$,  $\int_{\mathbb{R}}B_2(v)^2 dv$ and  $\int_{\mathbb{R}}C_2(v)^2 dv$.\\
We obtain estimates by replacing $t$ by $0$ in the results of 1.\\
\\
3.  Estimation of  $\int_{\frac{-s}{2}}^{-d}\Delta k^s_t(v)^2 dv$,  $\int_{\mathbb{R}}A_{3}(v)^2 dv$,  $\int_{\mathbb{R}}B_3(v)^2 dv$, 
$\int_{\mathbb{R}}C_3(v)^2 dv$, $\int_{\mathbb{R}}D(v)^2 dv$ and  $\int_{\mathbb{R}}E(v)^2 dv$.\\
 Denote $$G_4(z)\ :=\  z^{H-K}\Bigl(\hat{F}(-z)\ -\ 1\Bigr).$$ For $u\in(\frac{-s}{2},-d)$, we have that 
$\Delta k^s_t(u)\ =\ (u+s)^{H-K}\left(G_4\bigl(\frac{t-u}{u+s}\bigr)\ -\ G_4\bigl(\frac{-u}{u+s}\bigr)\right).$\\
\\
From (\ref{e:2.20}) and (\ref{e:2.181}), it follows that $\frac{dG_4}{dz}(z) \ =\ (H-K)z^{H-K-1}\left((1+z)^{H+K-1}-1\right).$\\
\\
Hence, by the mean value theorem, there exist $\theta_{s,u} \in \bigl(\frac{-u}{u+s},\frac{t-u}{u+s}\bigr)\subset \left(0,\frac{t+\frac{s}{2}}{\frac{s}{2}}\right)$ and $\eta_{s,u}\in(0,\theta_{s,u})$, such that 
\begin{eqnarray*}&&\Delta k^s_t(u)\ =\ (u+s)^{H-K}\left(\frac{t}{u+s}\right)(H-K)\theta_{s,u}^{H-K-1}\bigl((1+\theta_{s,u})^{H+K-1}-1\bigr)\\
&&=\ (u+s)^{H-K}\left(\frac{t}{u+s}\right) (H-K)\theta_{s,u}^{H-K-1}(H+K-1) \theta_{s,u}(1+\eta_{s,u})^{H+K-2}.\end{eqnarray*} Thus, \begin{equation}|\Delta k^s_t(u)| \ \leq\ |H-K||H+K-1|(u+s)^{H-K-1}t \left(\frac{-u}{u+s}\right)^{H-K-1}\frac{t-u}{u+s}.\label{proof3eq1}\end{equation}
In particular, we have that
 \begin{eqnarray}|\Delta k^s_t(u)| 
&\leq& |H-K||H+K-1|\left(\frac{s}{2}\right)^{-1}t \frac{t+d}{d}(-u)^{H-K}.\label{e:6.81}\end{eqnarray} 
So, for $K=\frac{1}{2}$, we obtain that 
 $\int^{-d}_{\frac{-s}{2}}\Delta k^s_t(u)^2 du\ \leq\ \frac{t^2}{8H}(\frac{t+d}{d})^2 s^{2H-2}.$\\
\\
Also, from (\ref{proof3eq1}), it follows for $u\in\bigl(\frac{-s}{2},-d\bigr)$, that \begin{equation}|\Delta k^s_t(u)| \ \leq\ |H-K||H+K-1|t\frac{t+\frac{s}{2}}{\frac{s}{2}}(-u)^{H-K-1}.\label{e:6.8}\end{equation}
Denote $$G_5(z)\ :=\ F\left(\frac{3}{2}-K,H-K,H-K+1,z\right).$$
Let $K\neq\frac{1}{2}$. For $H>K$ and $s>\max(2t,4d)$, we obtain by using (\ref{e:6.8}) and (\ref{e:2.263}), that\\
\\
$\int_{\mathbb{R}}A_3(v)^2dv\ \leq\  (H-K)^2 t^2 \left(\frac{t+\frac{s}{2}}{\frac{s}{2}}\right)^2     \int_{-\infty}^{\frac{-s}{2}}\left(\int^{-d}_{\frac{-s}{2}}(-d-u)^{H-K-1}(u-v)^{K-\frac{3}{2}}du\right)^2 dv$\\
\\ 
$\leq\ t^2 \left(\frac{t+\frac{s}{2}}{\frac{s}{2}}\right)^2   G^{\ast 2}_5   \frac{(\frac{s}{2}-d)^{2H-2}}{2-2K} \ \leq\ \frac{32 G^{\ast 2}_5 }{1-K}t^2 s^{2H-2}$.\\
\\
Similarly, for $H< K$, we have by using (\ref{e:6.81}) and (\ref{e:2.263}), that\\
\\
$\int_{\mathbb{R}}A_3(v)^2dv\ \leq\ \frac{2 G^{\ast 2}_6 t^2}{(H-K+1)^2 (1-K)}(\frac{t+d}{d})^2 s^{2H-2}$, where $$G_6(z)\ :=\ F\left(\frac{3}{2}-K,H-K+1,H-K+2,z\right).$$
Let $$G_7(z)\ :=\ F\bigl(2(K-H+1),1,2K+1,z\bigr).$$
For $K>\frac{1}{2}$  and $s>2t$, we have by using (\ref{e:6.8}) and then (\ref{e:2.264}) twice, that\\
\\
$\int_{\mathbb{R}}B_3(v)^2dv\ \leq\ t^2 \left(\frac{t+\frac{s}{2}}{\frac{s}{2}}\right)^2  \int^{-d}_{\frac{-s}{2}}\left( \int^{-d}_v (-u)^{H-K-1}(u-v)^{K-\frac{3}{2}}du  \right)^2 dv\ \leq$\\
\\
$ t^2 \left(\frac{t+\frac{s}{2}}{\frac{s}{2}}\right)^2   \frac{G^{\ast 2}_{1}}{\left(K-\frac{1}{2}\right)^2} \frac{1}{2K}\left(\frac{s}{2}-d\right)^{2K}\left(\frac{s}{2}\right)^{2(H-K-1)}G_7^{\ast}\ \leq\ \frac{8 G^{\ast 2}_1 G^{\ast}_7}{\left(K-\frac{1}{2}\right)^2 K}t^2 s^{2H-2}$,\\
\\
where $G_1$ is defined as in (\ref{ucdc}).\\
\\
Let $K<\frac{1}{2}$. It holds that \begin{eqnarray*}\frac{d^2 G_4}{d^2 z}(z) &=&  (H-K-1)(H-K)z^{H-K-2}\left((1+z)^{H+K-1}-1\right)\\
&& +\ (H-K)(H+K-1)z^{H-K-1}(1+z)^{H+K-2}.\end{eqnarray*}
For $u\in(\frac{-s}{2},-d)$, there exists $\theta_{s,u} \in \bigl(\frac{-u}{u+s},\frac{t-u}{u+s}\bigr)\subset \left(0,\frac{t+\frac{s}{2}}{\frac{s}{2}}\right)$, such that \\
\\
$\big|\frac{-d}{du}\Delta k^s_t (u)\big|\ \leq\ |H-K|(u+s)^{H-K-1}\Big|G_4\bigl(\frac{t-u}{u+s}\bigr)\ -\ G_4\bigl(\frac{-u}{u+s}\bigr)\Big|\ +$\\
\\
$ (u+s)^{H-K-2}t \Big|\frac{dG_4}{dz}\bigl(\frac{t-u}{u+s}\bigr)\Big|\ +\ (u+s)^{H-K-2}s \frac{t}{u+s}\Big|\frac{d^2 G_4}{d^2 z}(\theta_{s,v})\Big| \ \leq$\\
\\
$(H-K)^2|H+K-1|\frac{t+\frac{s}{2}}{\frac{s}{2}}\bigl(\frac{s}{2}\bigr)^{-1}t (-u)^{H-K-1}     \ +\  |H-K||H+K-1|\frac{t+\frac{s}{2}}{\frac{s}{2}}\bigl(\frac{s}{2}\bigr)^{-1} t (t-u)^{H-K-1}\ +$\\
\\
$|H-K-1||H-K||H+K-1|\frac{t+\frac{s}{2}}{\frac{s}{2}}s t \bigl(\frac{s}{2}\bigr)^{-1}(-u)^{H-K-2}\ +$\\
\\
$|H-K||H+K-1|s^{H+K-1} \bigl(\frac{s}{2}\bigr)^{-H-K}t(-u)^{H-K-1}.$\\
\\
Let $$G_8(z)\ :=\ F\left(1+K-H,K+\frac{1}{2},K+\frac{3}{2},z\right),$$ $$G_9(z)\ :=\ F\bigl(2(K+1-H),1,2K+3,z\bigr),$$  $$G_{10}(z)\ :=\ F\left(2+K-H,K+\frac{1}{2},K+\frac{3}{2},z\right)$$ and $$G_{11}(z)\ := \ F\bigl(2(K+2-H),1,2K+3,z\bigr).$$  
First, for $\hat{t}\in\{0,t\}$, we have by using (\ref{e:2.264}) twice, that\\
\\
$\int^{-d}_{\frac{-s}{2}}\left(\int^{-d}_v (\hat{t}-u)^{H-K-1}(u-v)^{K-\frac{1}{2}}du\right)^2 dv \ \leq\ \left(\frac{G^{\ast}_8}{K+\frac{1}{2}}\right)^2 \frac{G_9^{\ast}}{2K+2}\bigl(\frac{s}{2}\bigr)^{2H}.$\\
\\
Second, in the same way, we obtain that\\
\\
$\int^{-d}_{\frac{-s}{2}}\left(\int^{-d}_v (-u)^{H-K-2}(u-v)^{K-\frac{1}{2}}du\right)^2 dv \ \leq\ 
\left(\frac{G^{\ast}_{10}}{K+\frac{1}{2}}\right)^2 \frac{G^{\ast}_{11}}{2K+2}\bigl(\frac{s}{2}\bigr)^{2H-2}.$\\
\\
By using (\ref{sumequation}) with $n=4$, we obtain for $s>2t$, that \\
\\
$\int_{\mathbb{R}}C_3(v)^2dv \ \leq\  \frac{4}{\left(K-\frac{1}{2}\right)^2}\left(\frac{t+\frac{s}{2}}{\frac{s}{2}}\right)^2 t^2   \bigl(\frac{s}{2}\bigr)^{-2}   \left(\frac{G^{\ast}_8}{K+\frac{1}{2}}\right)^2 \frac{G^{\ast}_9}{2K+2}\bigl(\frac{s}{2}\bigr)^{2H}\ +$\\
\\
$\frac{4}{\left(K-\frac{1}{2}\right)^2}\left(\frac{t+\frac{s}{2}}{\frac{s}{2}}\right)^2 t^2 \bigl(\frac{s}{2}\bigr)^{-2} \frac{G^{\ast 2}_8}{\left(K+\frac{1}{2}\right)^2} \frac{G^{\ast}_9}{2K+2}\bigl(\frac{s}{2}\bigr)^{2H}\ +$\\
\\
$\frac{4}{\left(K-\frac{1}{2}\right)^2}  (H-K-1)^2 \left(\frac{t+\frac{s}{2}}{\frac{s}{2}}\right)^2 s^2 t^2  \bigl(\frac{s}{2}\bigr)^{-2}    \left(\frac{G^{\ast}_{10}}{K+\frac{1}{2}}\right)^2 \frac{G^{\ast}_{11}}{2K+2}(\frac{s}{2})^{2H-2}\ +$\\
\\
$\frac{4}{\left(K-\frac{1}{2}\right)^2}      s^{2(H+K-1)} t^2 \bigl(\frac{s}{2}\bigr)^{-2K-2H}       \left(\frac{G^{\ast}_8}{K+\frac{1}{2}}\right)^2 \frac{G^{\ast}_9}{2K+2}\bigl(\frac{s}{2}\bigr)^{2H}\ \leq\ \frac{\left(68 G^{\ast 2}_8 G^{\ast}_9 +512 G^{\ast 2}_{10}G^{\ast}_{11}\right)}{\left(K-\frac{1}{2}\right)^2 \left(K+\frac{1}{2}\right)^2}t^{2}s^{2H-2}.$\\
\\
By using (\ref{e:6.81}), we have that $|\Delta k^s_t(-d)|\ \leq\ |H-K||H+K-1|(\frac{s}{2})^{-1}t\frac{t+d}{d}d^{H-K}$. \\
 \\
Hence, $\int_{\mathbb{R}}D(v)^2 dv\ \leq\ \frac{2 d^{2(H-K)}}{K(K-\frac{1}{2})^2}(\frac{t+d}{d})^2 t^2 s^{2K-2}$.\\
\\
From (\ref{e:6.8}), it follows that $|\Delta k^s_t(\frac{-s}{2})|\ \leq\ |H-K||H+K-1|t\frac{t+\frac{s}{2}}{\frac{s}{2}}(\frac{s}{2})^{H-K-1}$.
\\
Hence for $s>2t$, it holds that $\int_{\mathbb{R}}E(v)^2 dv\ \leq\ \frac{8}{K(K-\frac{1}{2})^2}t^2 s^{2H-2}$.\\
\\
4. Estimation of  $\int^{\frac{-s}{2}}_{-s}\Delta g^s_t(v)^2 dv$,  $\int_{\mathbb{R}}A_{4}(v)^2 dv$,  $\int_{\mathbb{R}}B_4(v)^2 dv$ and  $\int_{\mathbb{R}}C_4(v)^2 dv$.\\
 By using (\ref{e:2.18}), we have that \begin{eqnarray*}G_{12}(z) &:=& z^{H-K}F(2H,H-K,H-K+1,z)\\
&=& \Bigl(\frac{1-z}{z}\Bigr)^{K-H}\hat{F}\Bigl(\frac{z}{z-1}\Bigr).\end{eqnarray*} From (\ref{e:2.20}) and (\ref{e:2.181}), it follows that \begin{eqnarray*}\frac{dG_{12}}{dz}(z) &=&  (H-K)z^{H-K-1}(1-z)^{-2H}.\end{eqnarray*}
Let  $u\in\bigl(-s,\frac{-s}{2}\bigr)$. There exists  $\theta_{s,u}\in\bigl(\frac{-u}{s},\frac{t-u}{s+t}\bigr)\ \subset\ \bigl(\frac{1}{2},1\bigr)$, such that\begin{eqnarray*}
\Delta g^s_t(u)&=& (u+s)^{H-K}\left(G_{12}\left(\frac{t-u}{s+t}\right)\ -\ G_{12}\left(\frac{-u}{s}\right)\right)\\
&=& (u+s)^{H-K}\left(\frac{(s+u)t}{(s+t)s}\right)(H-K)\theta_{s,u}^{H-K-1}(1-\theta_{s,u})^{-2H}.\end{eqnarray*}
In particular, $|\Delta g^s_t(u)|\ \leq\ \frac{|H-K|t}{(s+t)^{1-2H}s}(\frac{1}{2})^{H-K-1}(u+s)^{1-K-H}$.\\
\\
Hence, for $K=\frac{1}{2}$ and $s>t$, we have that $\int^{\frac{-s}{2}}_{-s}\Delta g_t^s(u)^2 du\leq\frac{1}{1-H}t^2 s^{2H-2}$.\\
\\
Denote $$G_{13}(z)\ :=\ F\left(\frac{3}{2}-K,1,3-K-H,z\right).$$ For  $K\neq\frac{1}{2}$ and $s>t$, (\ref{e:2.263}) yields that  $\int_{\mathbb{R}}A_4(v)^2dv\ \leq\ \frac{2 G^{\ast 2}_{13}t^2}{(1-K)(2-K-H)^2}s^{2H-2}$.\\
\\
Similarly, for $K>\frac{1}{2}$ and $s>t$, we obtain that $\int_{\mathbb{R}}B_4(v)^2dv\ \leq\ \frac{2 G^{\ast 2}_{14}}{\left(K-\frac{1}{2}\right)^2 K}t^2 s^{2H-2}$, where $$G_{14}(z) \ :=\ F\left(H+K-1,1,K+\frac{1}{2},z\right).$$
Furthermore, we  have that $$\frac{d^2 G_{12}}{d^2 z}(z) =
(H-K)\left((H-K-1)z^{H-K-2}(1-z)^{-2H}+ 2H z^{H-K-1}(1-z)^{-2H-1}\right).$$
So for  $u\in(-s,\frac{-s}{2})$, it holds that\\
\\
$\Big|\frac{d}{du}\Delta g^s_t(u)\Big| \leq\ |H-K|(u+s)^{H-K-1}\big|G_{12}(\frac{t-u}{s+t})\ -\ G_{12}(\frac{-u}{s})\big|\ +$\\
\\
$(u+s)^{H-K}\big|\frac{dG_{12}}{dz}(\frac{t-u}{s+t})(\frac{-1}{s+t})-\frac{dG_{12}}{dz}(\frac{-u}{s})\frac{-1}{s}\big|\ \leq$\\
\\
$(H-K)^2 \frac{t (\frac{1}{2})^{H-K-1}}{(s+t)^{1-2H}s}(u+s)^{-H-K}\ +\ (u+s)^{H-K}\bigl|\frac{dG_{12}}{dz}(\frac{t-u}{s+t})-\frac{dG_{12}}{dz}(\frac{-u}{s})\bigr|\frac{1}{s+t}$\\
\\
$\ +\ (u+s)^{H-K}\big|\frac{dG_{12}}{dz}\bigl(\frac{-u}{s}\bigr)\big|\frac{t}{(s+t)s}\ \leq$\\
\\
$\frac{(H-K)^2 t(\frac{1}{2})^{H-K-1}}{(s+t)^{1-2H}s}(u+s)^{-H-K}\  +\ \frac{|H-K||H-K-1|t(\frac{1}{2})^{H-K-2}}{s(s+t)^{2-2H}}(u+s)^{1-H-K}\ + $\\
\\
$ \frac{|H-K|2H t(\frac{1}{2})^{H-K-1}}{s(s+t)^{1-2H}}(u+s)^{-H-K}\ + \ \frac{|H-K|t(\frac{1}{2})^{H-K-1}}{(s+t)s^{1-2H}}(u+s)^{-H-K}\ =$\\
\\
$\left(\frac{(H-K)^2t(\frac{1}{2})^{H-K-1} }{(s+t)^{1-2H}s}\ +\  \frac{|H-K|2H t(\frac{1}{2})^{H-K-1}}{s(s+t)^{1-2H}}\ +\ \frac{|H-K|t(\frac{1}{2})^{H-K-1}}{(s+t)s^{1-2H}}\right)(u+s)^{-H-K}\ +$\\
\\
$\frac{|H-K||H-K-1|t(\frac{1}{2})^{H-K-2}}{s(s+t)^{2-2H}}(u+s)^{1-H-K}$.\\
\\
Denote $$G_{15}(z)\ :=\ F\left(K+H,1,K+\frac{3}{2},z\right)$$ and  $$G_{16}(z)\ :=\ F\left(K+H-1,1,K+\frac{3}{2},z\right).$$
Then, for $K<\frac{1}{2}$ and $s>t$, we obtain by using (\ref{e:2.263}) and (\ref{sumequation}) with $n=3$, that\\
\\
$\int_{\mathbb{R}}C_4(v)^2dv\ \leq\ \left(\frac{1}{K-\frac{1}{2}}\right)^2 \left(\frac{t\left(\frac{1}{2}\right)^{H-K-1} }{(s+t)^{1-2H}s}\ +\  \frac{2 t\left(\frac{1}{2}\right)^{H-K-1}}{s(s+t)^{1-2H}}\ +\ \frac{t\left(\frac{1}{2}\right)^{H-K-1}}{(s+t)s^{1-2H}}\right)^2 \frac{G^{\ast 2}_{15}}{\left(K+\frac{1}{2}\right)^2} \frac{\left(\frac{s}{2}\right)^{2-2H}}{2K+2}\ +$\\
\\
$\left(\frac{1}{K-\frac{1}{2}}\right)^2 \left( \frac{|H-K-1|t\left(\frac{1}{2}\right)^{H-K-2}}{s(s+t)^{2-2H}}\right)^2 \frac{G^{\ast 2}_{16}}{\left(K+\frac{1}{2}\right)^2}\frac{\left(\frac{s}{2}\right)^{4-2H}}{2K+2}\ \leq\ \frac{\left(63 G^{\ast 2}_{15}+ 4 G^{\ast 2}_{16}\right)}{\left(K-\frac{1}{2}\right)^2 \left(K+\frac{1}{2}\right)^2}t^2 s^{2H-2}$.\\
\\
5. Estimation of  $\int^{\frac{-s}{2}}_{-s}\Delta h^s_t(v)^2 dv$,  $\int_{\mathbb{R}}A_{5}(v)^2 dv$,  $\int_{\mathbb{R}}B_5(v)^2 dv$ and  $\int_{\mathbb{R}}C_5(v)^2 dv$.
 
For $u\in\bigl(-s,\frac{-s}{2}\bigr)$, we have that $|\Delta h^s_t(u)|\  \leq\ |H-K|t(\frac{s}{2})^{H-K-1}$.\\
For $K=\frac{1}{2}$, it follows that $\int^{\frac{-s}{2}}_{-s}\Delta h^s_t(u)^2 du\ \leq t^2 s^{2H-2}$.\\
\\
Let $$G_{17}(z)\ :=\  F\left(\frac{3}{2}-K,1,2,z\right).$$
For $K\neq\frac{1}{2}$, we obtain by using  (\ref{e:2.263}), that $\int_{\mathbb{R}}A_5(v)^2 dv\ \leq\ \frac{2 G^{\ast 2}_{17}}{1-K}t^2 s^{2H-2}.$\\
\\
Moreover, for $K>\frac{1}{2}$, we have that 
$\int_{\mathbb{R}}B_5(v)^2 dv\ \leq\ \frac{2 t^2}{K\left(K-\frac{1}{2}\right)^2}s^{2H-2}$.\\
\\
For $u\in\left(-s,\frac{-s}{2}\right)$, we have that $\big|\frac{d}{du}\Delta h^s_t(u)\big|\leq |H-K||H-K-1|t\bigl(\frac{s}{2}\bigr)^{H-K-2}$.\\
\\
For $K<\frac{1}{2}$, it follows that  $\int_{\mathbb{R}}C_5(v)^2 dv\ \leq \ \frac{8 \Gamma\left(K-\frac{1}{2}\right)^2 }{\Gamma\left(K+\frac{3}{2}\right)^2 }t^2 s^{2H-2}$.\\
\\
By combining (\ref{e:6.5}) and (\ref{e:6.6}), respectively, with these estimates, we obtain (\ref{e:6.21}) where $$c_2\left(\frac{1}{2},H,t,d\right)\: =\ 2\cdot ^{\ast}\!G_2^2 (t+d)^{2H+2}\ +\ \left(\frac{(t+d)^2}{4Hd^2}+\frac{2}{1-H}+2\right)t^2$$ and\\
\\
$c_2(K,H,t,d)\ =\ \left(\frac{20\cdot ^{\ast}\!G_2^2 G^{\ast 2}_3}{\Gamma\left(K-\frac{1}{2}\right)^2(1+H-K)^2(1-K)}+\frac{20\cdot ^{\ast}\!G_2^2 \Gamma(H-K+1)^2}{\Gamma\left(H+\frac{3}{2}\right)^2}\right)(t+d)^{2H+2}\ +$\\
\\
\hbox{}\qquad\qquad\qquad$\left(\frac{\max\left(160G^{\ast 2}_5,\frac{10G^{\ast 2}_6 (t+d)^2}{(H-K+1)^2 d^2}\right)+10 G^{\ast 2}_{17}+\frac{10 G^{\ast 2}_{13}}{(2-K-H)^2}}{(1-K)\Gamma\left(K-\frac{1}{2}\right)^2}\  +\
\frac{10 G^{\ast 2}_{14}+40G^{\ast 2}_1G^{\ast}_7+10}{K\Gamma\left(K+\frac{1}{2}\right)^2}\right)t^2$,\quad $K>\frac{1}{2}$.\\
\\
In the same way, by combing (\ref{e:6.7}) with the estimates, we obtain (\ref{e:6.22}) with \\
\\
$c_3(K,H,t,d)\ =\ \left(\frac{24\cdot ^{\ast}\!G^{2}_2 G^{\ast 2}_3}{\Gamma\left(K-\frac{1}{2}\right)^2(1+H-K)^2(1-K)}+\frac{816\cdot ^{\ast}\!G^2_2 \Gamma(H-K+1)^2}{(1+H-K)^2\Gamma\left(H+\frac{3}{2}\right)^2}\right)(t+d)^{2H+2}\ +$\\
\\
\hbox{}\quad\qquad\qquad\qquad$\left(\frac{\max\left(192 G^{\ast 2}_{5},\frac{12 G^{\ast 2}_{6}(t+d)^{2}}{(H-K+1)^2 d^2} \right)\ +\ \frac{12 G^{\ast 2}_{13}}{(2-H-K)^2}+12G^{\ast 2}_{17}}{\Gamma\left(K-\frac{1}{2}\right)^2(1-K)}\right)t^2\ +\ \frac{48}{K\Gamma\left(K+\frac{1}{2}\right)^2}t^2\ +$\\
\\
\hbox{}\quad\qquad\qquad\qquad
$\left(\frac{408 G^{\ast 2}_8 G^{\ast}_9 + 3072 G^{\ast 2}_{10}G^{\ast}_{11}+378 G^{\ast 2}_{15}+24 G^{\ast 2}_{16}+48}{\Gamma\left(K+\frac{3}{2}\right)^2}\right)t^2$\\
\\
and  $$c_4(K,H,t,d)\ =\
  \frac{12 d^{2(H-K-1)} (t+d)^2 }{\Gamma\left(K+\frac{1}{2}\right)^2 K}t^2.$$
Eventually, by combining (\ref{e:6.3}) and (\ref{e:6.4})  with the results in \textbf{1.} and \textbf{2.}, we obtain the following: For $K=\frac{1}{2}$, we have that \begin{eqnarray*}\left(\frac{\Gamma\bigl(H+\frac{1}{2}\bigr)}{C(H)}\right)^2 E\bigl[Z^{H,s}_t- Z^{H}_t\bigr]^2& \leq &\left(c_1\left(\frac{1}{2},H,t\right)\ +\ c_2\left(\frac{1}{2},H,t,d\right)\right)s^{2H-2},\\
& &\qquad\qquad  s>2t+4d+1.\end{eqnarray*}
For $K>\frac{1}{2}$, we have that  \begin{eqnarray*}\frac{\Gamma(H-K+1)^2}{2C(H)^2} E\bigl[Z^{H,s}_t-Z^{H}_t\bigr]^2 &\leq&\bigl(c_1(K,H,t)+c_2(K,H,t,d)\bigr)s^{2H-2},\\ & &\qquad\qquad s>2t+4d+1.\end{eqnarray*}
For $K<\frac{1}{2}$, it holds that \begin{eqnarray*}\frac{\Gamma(H-K+1)^2}{2C(H)^2}\cdot  E\bigl[Z^{H,s}_{t}-Z^{H}_t\bigr]^2 &\leq& \bigl(c_1(K,H,t)\ +\ c_3(K,H,t,d)\bigr)s^{2H-2}\ +\\
 & &c_4(K,H,t,d)s^{2K-2},\ s>2t+4d+1.\end{eqnarray*}
\end{proof}
The (generalized) Mandelbrot-Van\,Ness representation is a direct consequence of Theorem \ref{theorem1}, and hence of the (generalized) Molchan-Golosov representation:
\begin{cor}For every $K\in(0,1)$, the process $\bigl(Z^H_t\bigr)_{t\in\mathbb{R}}$ is an $H$-fBm.\end{cor}
\begin{proof}
For $t\in\mathbb{R}$, it follows from Theorem \ref{theorem1}, that $\lim_{s\to\infty}E\bigl[Z^{H,s}_t - Z^H_t\bigr]^2=0$. Hence, for all $t,t'\in\mathbb{R}$, we have that $E\bigl[Z^{H}_t\cdot Z^H_{t'}\bigr]= \lim_{s\to\infty}E\bigl[Z^{H,s}_t\cdot Z^{H,s}_{t'}\bigr] = \frac{1}{2}\left(|t|^{2H}+|t'|^{2H}-|t-t'|^{2H}\right)$.  \end{proof}
\textbf{Acknowledgements.} Thanks to my supervisor Esko Valkeila for suggesting the topic stu\--died here and for useful remarks. Thanks are also due to the referee of Stochastic Processes and their Applications for a careful reading in connection with \cite{Jo}. I would like to thank DYNSTOCH Network and the Finnish Graduate School in Stochastics (FGSS) for financial support. 


\begin{thebibliography}{100}
\bibitem{Ab}Abramowitz, M., Stegun, I.A., Handbook of Mathematical Functions, Dover Publications Inc., New York, 1965.
\bibitem{De2}Decreusefond, L., \"Ust\"unel, A.S., Stochastic Analysis of the Fractional Brownian Motion, Potential Analysis 10 (1999) 177-214.
\bibitem{Er}Erd\'elyi, A., Magnus, W., Oberhettinger, F., Tricomi, F.G., Higher Transcendental Functions. Volume 1, McGraw-Hill Book Company, Inc., 1953.
\bibitem{Jo}Jost, C., Transformation formulas for fractional Brownian motion, Stochastic Processes and their Applications 116 (2006) 1341-1357.
\bibitem{Ma}Mandelbrot, B.B., Van\,Ness, J.W., Fractional Brownian Motions, Fractional Noises and Applications, SIAM Review 10(4) (1968) 422-437.
\bibitem{Mo1}Molchan, G., Golosov, J., Gaussian Stationary Processes with Asymptotic Power Spectrum, Soviet. Math. Doklady 10(1) (1969) 134-137.
\bibitem{Pi1}Pipiras, V., Taqqu, M.S., Fractional calculus and its connections to fractional Brown\-ian motion, in: P. Doukhan, G. Oppenheim, and M.S. Taqqu (Eds.), Theory and Applications of Long-Range Dependence, Birkh\"auser, 2003, pp.165-201.
\bibitem{Pi2}Pipiras, V., Taqqu, M.S., Deconvolution of fractional Brownian motion, Journal of Time Series Analysis 23 (2002) 487-501.
\bibitem{Pi4}Pipiras, V., Taqqu, M.S., Integration questions related to fractional Brownian motion, Probability Theory and Related Fields 118(2) (2000) 251-291. \bibitem{Sa}Samko, S.G., Kilbas, A.A., Marichev, O.I., Fractional Integrals and Derivatives: Theory and Applications, Gordon and Breach Science Publishers, London, 1993.
\bibitem{Sam}Samorodnitsky, G., Taqqu, M.S., Stable non-Gaussian random processes: stochastic models with infinite variance, Chapman and Hall, New York, 1994.\end{thebibliography}
\end{document}